\documentclass[11pt, a4paper, reqno]{amsart}
\usepackage{amssymb, array, amsmath, amscd, pdfpages, enumerate, amsthm, fixltx2e, setspace, hyperref,mathtools}

\usepackage[utf8]{inputenc}
\usepackage{amssymb}

\usepackage{graphicx, caption} 
\usepackage{tikz}
\usepackage{bm}

\newcommand*{\R}{{\mathbb{R}}}

\newcommand*{\N}{{\mathbb{N}}}

\newcommand*{\Abs}[2][default]{\ifthenelse{\equal{#1}{default}}{\left\lvert#2\right\rvert}{\ldelim{#1}{\lvert}#2\rdelim{#1}{\rvert}}}
\newcommand*{\Norm}[2][default]{\ifthenelse{\equal{#1}{default}}{\left\lVert#2\right\rVert}{\ldelim{#1}{\lVert}#2\rdelim{#1}{\rVert}}}

\newcommand*{\Iprod}[3][default]{\ifthenelse{\equal{#1}{default}}{\left\langle#2,#3\right\rangle}{\ldelim{#1}{\langle}#2,#3\rdelim{#1}{\rangle}}}
\newcommand*{\Dualpair}[3][default]{\ifthenelse{\equal{#1}{default}}{\left\langle#2,#3\right\rangle}{\ldelim{#1}{\langle}#2,#3\rdelim{#1}{\rangle}}}

\newcommand*{\ddb}[2][1]{\ifthenelse{\equal{#1}{1}}{\frac{d}{d#2}}{\frac{d^{#1}}{d#2^{#1}}}}
\newcommand*{\pd}[3][1]{\ifthenelse{\equal{#1}{1}}{\frac{\partial{#2}}{\partial{#3}}}{\frac{\partial^{#1}{#2}}{\partial#3^{#1}}}}

\newcommand{\blue}[1]{{\color{blue} #1}}
\renewcommand{\blue}[1]{#1}

\newcommand*\lenv{{\hbox{\raisebox{-.15ex}{\rotatebox[origin=c]{50}{$\smallsmile$}}\kern-8.65pt\rotatebox[origin=c]{-25}{$\smallsetminus$}}}}
\newcommand*\uenv{{\hbox{\raisebox{-.0ex}{\rotatebox[origin=c]{-45}{$\smallfrown$}}\kern-5.6pt\raisebox{.2ex}{\rotatebox[origin=c]{-5}{\scriptsize\slash}}}}\,\kern+1.5pt}
\newcommand{\sleq}{\preccurlyeq}
\newcommand{\sgeq}{\succcurlyeq}
\newcommand*{\Gg}{{\mathcal{G}}} 
\newcommand*{\Ss}{{\mathcal{S}}}

\newcommand*{\Aa}{{\mathcal{A}}}
\newcommand*{\Mm}{{\mathcal{M}}}

\newcommand*{\lupp}[1]{\kern0pt^{\kern0pt u \kern0pt}\kern1pt#1}
\newcommand*{\llow}[1]{\kern0pt^{\kern0pt l \kern0pt}\kern1pt#1}
\newcommand*{\rupp}[1]{#1\kern1pt^{\kern0pt u \kern0pt}\kern0pt}
\newcommand*{\rlow}[1]{#1\kern1pt^{\kern0pt l \kern0pt}\kern0pt}
\newcommand*{\ul}[1]{\kern0pt^{\kern0pt u \kern0pt}\kern0pt#1\kern1pt^{\kern0pt l \kern0pt}\kern0pt}
\newcommand*{\lu}[1]{\kern0pt^{\kern0pt l \kern0pt}\kern0pt#1\kern0pt^{\kern0pt u \kern0pt}\kern0pt}
\newcommand*{\bproofname}{Proof}
\newenvironment{bproof}[1][\bproofname]{\begin{proof}[#1]}{\end{proof}}

\newtheorem{thm}{Theorem}[section]
\numberwithin{thm}{section}

\newtheorem{lemma2}[thm]{Lemma}

\theoremstyle{definition}
\newtheorem{defn}[thm]{Definition}
\newtheorem{remark2}[thm]{Remark}
\newtheorem{example2}[thm]{Example}
\numberwithin{equation}{section}

\begin{document}

\title[General mixed lattices]{General mixed lattices}

\thispagestyle{plain}

\author{Jani Jokela}
\address[J. Jokela]{Mathematics, Faculty of Information Technology and Communication Sciences, 
Tampere University, PO. Box 692, 33101 Tampere, Finland}
\email{jani.jokela@tuni.fi}


\begin{abstract}
A mixed lattice is a lattice-type structure consisting of a set with two partial orderings, and generalizing the notion of a lattice. 
Mixed lattice theory has previously been studied in various algebraic structures, such as groups and semigroups, while the more general notion of a mixed lattice remains unexplored. 
In this paper, we study the fundamental properties of mixed lattices and the relationships between the various properties, such as the one-sided associative, distributive and modular laws. We also give an alternative definition of mixed lattices and mixed lattice groups as non-commutative and non-associative algebras satisfying a certain set of postulates. The algebraic and the order-theoretic definitions are then shown to be equivalent. 
\end{abstract}

\maketitle

\section{Introduction}\label{sec:sec1}

We first give a short overview of the subject. All the terminology and concepts mentioned in the introduction will be given a more detailed explanation in Section \ref{sec:sec2}. 

The development of the theory of algebraic mixed lattice structures started with the work of Arsove and Leutwiler on mixed lattice semigroups \cite{ars3}, \cite{ars4}, \cite{ars}, and was later continued by Eriksson-Bique on mixed lattice groups \cite{eri1}, \cite{eri}. More recently, a systematic study of mixed lattice groups and vector spaces has been undertaken in \cite{jj1}, \cite{jj2}, \cite{jj3}. 
\blue{Algebraic mixed lattice structures appear in many branches of mathematics. 
Historically, mixed lattice structures arose from potential theory, where the set of superharmonic functions defined on some region in $\R^n$ forms a mixed lattice semigroup. 
Some other important examples of mixed lattice structures include the set of functions of bounded variation, and the set of functions that can be written as a difference of two concave positive functions \cite{ars}. 
Cone projections play an important role in convex analysis and its applications to optimization problems, and the author has shown in \cite{jj4} that the cone projection problem can be given an abstract formulation based on a mixed lattice structure on a vector space.  
Moreover, in functional analysis, mixed lattice vector spaces provide a natural setting for the study of asymmetric norms and the related topologies \cite{jj3}. In the present paper, we give additional examples from algebra and the theory of divisibility (see Examples \ref{esim1} and \ref{multiplicative}). Other examples of algebraic mixed lattice structures can be found in \cite{jj1}.} 

The theory developed in the articles cited above is a generalization of the theory of vector lattices and lattice ordered groups (for which we refer to \cite{lux} and \cite{birk}). This generalization is based on the idea of starting with a partially ordered structure (such as a group or a vector space) and introducing a second partial ordering, resulting in a more general ordered structure in which the two partial orderings are mixed.    
A similar generalization of lattice theory can be obtained by studying the mixed lattice order structure in its own right, without assuming any underlying algebraic structure. This will be the topic of the next section, where we study the basic properties of mixed lattices which are similar to the general properties of lattices. 
\blue{These include the one-sided associative and distributive laws.  
These properties are known in mixed lattice semigroups and groups but here we show that they are in fact equivalent (Theorem \ref{main1}). 
These results are of particular importance because a great deal of the theory of mixed lattice vector spaces depend on properties of Theorem \ref{main}(a) and Theorem \ref{main1}(a).  
For more on these developments, see \cite{jj2}.} 

A lattice can be defined either as a partially ordered set, or as an algebraic system with two binary operations that are commutative, idempotent and associative. In Section \ref{sec:sec3} we study a similar alternative definition of a mixed lattice as a non-commutative and non-associative algebra with two binary operations satisfying a set of axioms. We then show that the algebraic and the order-theoretic definitions are equivalent in a certain sense. We also give a similar axiomatic definition of a mixed lattice group. 
 
As this theory is still in its infancy, it goes without saying that there is a huge number of possibilities for further research. The purpose of the present paper is to serve as a starting point, and 
the author hopes that this introductory paper would 
inspire mathematicians to further develop the theory of mixed lattices.

\section{Mixed lattices}\label{sec:sec2}

We begin with a review of some elementary facts from lattice theory. A binary relation $\leq$ defined on a set $M$ is called a \emph{partial ordering} if the following properties hold for all $x,y,z\in M$.   
\[ 
\begin{array}{rcl}
x\leq x  & & \text{(Reflexivity)} \\
x\leq y \; \text{ and } \; y\leq x \; \text{ implies } \; y=x & & \text{(Antisymmetry)} \\
x\leq y \; \text{ and } \; y\leq z \; \text{ implies } \; x\leq z & & \text{(Transitivity)}  \\[2ex]
\end{array}
\]
A set $M$ together with a partial ordering $\leq$, denoted by $(M,\leq)$, is called a 
\emph{partially ordered set}. 
A partially ordered set $(M,\leq)$ is called a \emph{lattice} if $\sup\{x,y\}$ and $\inf \{x,y\}$ exist in $M$ for all $x,y\in M$. The standard notation for these is $\sup\{x,y\}=x\vee y$ and $\inf \{x,y\}=x\wedge y$. There are various identities and inequalities that hold in every lattice, and each of them have dual versions, obtained simply by reversing the ordering $\leq$ 
(the so-called \emph{duality principle}). We will collect some of these below. %
For more details on general lattice theory, including the proofs of the following results, we refer the reader to the standard references \cite{birk} and \cite{gra}. 

If $(M,\leq)$ is a lattice then the following hold for all $x,y,z\in M$.

\begin{equation}\label{lass}
(x\wedge y)\wedge z = x\wedge (y\wedge z) \; \textrm{ and } \; (x\vee y)\vee z = x\vee (y\vee z),
\end{equation}
\begin{equation}\label{ldist1}
(z\wedge x)\vee (z\wedge y) \leq z\wedge(x\vee y), 
\end{equation}
\begin{equation}\label{ldist2}
z\vee (x\wedge y) \leq (z\vee x)\wedge (z\vee y), 
\end{equation}
\begin{equation}\label{lmod1}
(z\wedge x)\vee (z\wedge y) \leq z\wedge((z\wedge x)\vee y), \\[2ex]
\end{equation}

The identities \eqref{lass} are called the \emph{associative laws}, the inequalities \eqref{ldist1} and \eqref{ldist2} are called the
the \emph{distributive inequalities}, and \eqref{lmod1} is called the \emph{modular inequality}. If equality holds in \eqref{ldist1} and \eqref{ldist2} for all $x,y,z\in M$, then the lattice $(M,\leq)$ is called \emph{distributive}. If equality holds in \eqref{lmod1} for all $x,y,z\in M$ then the lattice $(M,\leq)$ is called \emph{modular}. Moreover, distributivity implies modularity, but not conversely.

Suppose we have two partial orderings $\leq$ and $\sleq$ defined on a set $M$. 
Then for every $x,y\in M$ the \emph{mixed upper envelope} is defined by 
\begin{equation}\label{upperenv}
x\uenv y\,=\,\min \,\{\,w\in M: \; w\sgeq x \; \textrm{ and } \; w\geq y \,\},
\end{equation}
and the \emph{mixed lower envelope} is defined by  
\begin{equation}\label{lowerenv}
x\lenv y \,=\,\max \,\{\,w\in M: \; w\sleq x \; \textrm{ and } \; w\leq y \,\}, 
\end{equation}
where the minimum and maximum (if they exist) are taken with respect to the partial order $\leq$.

\begin{defn}\label{ml}
If $M$ is a set and  
$\leq$ and $\sleq$ are two partial orderings defined on $M$, then $\Mm=(M,\leq,\sleq)$ is called a \emph{mixed lattice} if the elements $x\uenv y$ and $x\lenv y$ defined by \eqref{upperenv} and \eqref{lowerenv} exist in $M$ for all $x,y\in M$. A subset $\Ss$ of $\Mm$ is called a \emph{mixed sublattice of} $\Mm$ if $x,y\in \Ss$ implies that $x\uenv y, x\lenv y \in \Ss$ (where the mixed envelopes are taken in $\Mm$), and the mixed envelopes of $\Ss$ are the restrictions of $\uenv$ and $\lenv$ to $\Ss$.
\end{defn}

We note that the operations $\uenv$ and $\lenv$ are not commutative, that is, $x\uenv y\neq y\uenv x$ and $x\lenv y\neq y\lenv x$, in general. 
However, if the partial orderings $\leq$ and $\sleq$ are identical then $x\uenv y=y\uenv x=\sup \{x,y\}$ and $x\lenv y=y\lenv x=\inf\{x,y\}$, and in this case, $(M,\leq,\sleq)=(M,\leq)$ is an ordinary lattice. Hence, a mixed lattice is a generalization of a lattice.

In Section \ref{sec:sec1}, we mentioned some of the algebraic mixed lattice structures. 
We recall that a commutative group $(G,*)$ together with a partial order $\leq$ is called a \emph{partially ordered group} if $\leq$ is such that for all $x,y\in G$
\begin{equation}\label{gorder}
x\leq y \, \implies x * z \leq y * z \quad \textrm{for all} \, z\in G. \\[1ex]
\end{equation}

We now have the following definition.

\begin{defn}\label{mlg}
Let $(G,*)$ be a commutative group equipped with two partial orderings $\leq$ and $\sleq$ such that $(G,*,\leq)$ and $(G,*,\sleq)$ are both partially ordered groups. If $(G,\leq, \sleq)$ is a mixed lattice, as defined in Definition \ref{ml}, then $\Gg=(G,*,\leq, \sleq)$ is called a \emph{mixed lattice group}.
\end{defn}

A \emph{mixed lattice vector space} is defined similarly as a partially ordered vector space $V$ with two partial orderings $\leq$ and $\sleq$ such that $(V,\leq,\sleq)$ is a mixed lattice. In this case, the group operation is the vector addition. Traditionally, the additive notation is used for the group operation in the theory of mixed lattice groups and semigroups, 
but there is no particular reason why the group operation should be addition (cf. Example \ref{multiplicative}), so we use here a generic symbol $*$ for the group operation. 
The definition of a mixed lattice semigroup is somewhat complicated, and we refer to \cite{ars} and \cite{jj1} for more details.


The following theorem lists the basic properties of mixed lattices that follow from the definition. These properties are known in mixed lattice groups and semigroups, but we will give here proofs that do not depend on any additional algebraic structure.

\begin{thm}\label{basic}
Let $\Mm=(M,\leq,\sleq)$ be a mixed lattice. The following hold.
\begin{enumerate}[\normalfont(a)]
\item
$x\uenv x=x$\, and \, $x\lenv x=x$\, for all $x\in\Mm$.
\item
$x\lenv y \sleq x\sleq x\uenv y$ \, and \, $x\lenv y \leq y\leq x\uenv y$ for all $x,y\in \Mm$.  
\item
$x\sleq u \; \textrm{ and } \; y\leq v \; \implies \; x\uenv y \leq u\uenv v \; \textrm{ and } \; x
\lenv y \leq  u\lenv v$. 
\item
$x\sleq y \, \iff \, x\uenv y = y \, \iff \, y\lenv x = x$.
\item
$(x\lenv y)\uenv x=x$ \, and \, $(x\uenv y)\lenv x=x$\, for all $x,y\in\Mm$. 
\end{enumerate}
\end{thm}

\begin{bproof}
(a) and (b) are evident by the definition of the mixed envelopes. 

(c)\;
Assume that $x\sleq u$  and  $y\leq v$ and define the sets 
$A=\{w\in M:w\sgeq x \; \textrm{and} \; w\geq y \}$ and 
$B=\{w\in M:w\sgeq u \; \textrm{and} \; w\geq v \}$. 
Then certainly $B\subseteq A$, and hence $\min A\leq \min B$, that is, $x\uenv y \leq u\uenv v$. A similar argument shows that $x\lenv y \leq u\lenv v$.

(d)\;
Let $x\sleq y$. Since $y\leq y$, it follows by parts (a), (b) and (c) that $y\leq x\uenv y\leq y\uenv y=y$. It follows that $y=x\uenv y$, by antisymmetry. Similarly, we have $x=x\lenv x\leq y\lenv x\leq x$, and so $y\lenv x=x$. The converse implications follow immediately from (b), since $x\sleq x\uenv y =y$ and $x=y\lenv x \sleq y$.

(e)\;
By (b) we have $x\lenv y\sleq x$, and so $(x\lenv y)\uenv x=x$ follows from (d). The other identity is proved similarly. 
\end{bproof}

The identities in Theorem \ref{basic}(e) are called the \emph{absorption laws}. 

Theorem \ref{basic}(d) gives a characterization for the partial order $\sleq$ in terms of the mixed envelopes. 
Next, we would like to have a similar characterization for the partial order $\leq$. Let us consider the following statement:

\begin{equation}\label{r0}
x\leq y \, \iff \, y\uenv x = y \, \iff \, x\lenv y = x. \\[2ex]
\end{equation}

If $y=y\uenv x$ then it follows immediately from Theorem \ref{basic}(b) that $x\leq y$, and similarly  $x=x\lenv y$ implies $x\leq y$. 
Unfortunately, it seems that we cannot get much further than this with the general Definition \ref{ml}  without assuming that the two partial orderings are connected in some way. 
First we give the following characterization of \eqref{r0} in terms of the two partial orderings. 

\begin{lemma2}\label{r0char}
The statement \eqref{r0} is equivalent to the following condition: If $y\sleq x$ and $x\leq y$ then $x=y$.
\end{lemma2}

\begin{bproof}
Assume first that \eqref{r0} holds. If $y\sleq x$ and $x\leq y$ then we have $y=x\lenv y$ and 
$x=x\lenv y$. Hence, $x=y$. Conversely, if the given condition holds and $x\leq y$ then by Theorem \ref{basic} we have $x=x\lenv x \leq x\lenv y \sleq x$, and this implies that $x=x\lenv y$. Similarly, $y\sleq y\uenv x\leq y\uenv y=y$ implies $y=y\uenv x$. %
\end{bproof}

A somewhat stronger additional requirement is the following, so-called \emph{pre-regularity condition}:

\begin{equation}\label{prereg}
x\sleq y  \; \implies \; x\leq y. \\[2ex]
\end{equation}

The condition \eqref{prereg} is sufficient for the equivalences in \eqref{r0}.

\begin{thm}\label{prereg_r0}
If the pre-regularity condition \eqref{prereg} holds in $\Mm=(M,\leq,\sleq)$ then the equivalences given in \eqref{r0} hold too. 
Moreover, if \eqref{r0} holds, then the absorption identities $x\uenv (y\lenv x)=x$ \, and \, $x\lenv (y\uenv x)=x$\, hold for all $x,y\in M$.
\end{thm}

\begin{bproof}
It was already noted above that the implications $y=y\uenv x$  $\implies$  $x\leq y$ \, and \, $x=x\lenv y$  $\implies$  $x\leq y$ 
hold (even without assuming the validity of \eqref{prereg}). Conversely, if the pre-regularity condition \eqref{prereg} holds then $y\sleq x$ and $x\leq y$ implies $y\leq x$ and $x\leq y$, and so $x=y$ by antisymmetry. Then \eqref{r0} holds, by Lemma \ref{r0char}. 
The absorption identities are then proved in the same way as in Theorem \ref{basic}(e).  
\end{bproof}

\begin{remark2}
The condition \eqref{r0} does not hold in general (see \cite[Example 2.25]{jj1}), and \eqref{r0} does not imply the pre-regularity condition. An example of the latter can be obtained from   
\cite[Example 4.10]{jj1} by exchanging the roles of the orders $\leq$ and $\sleq$.
\end{remark2}

We next consider several additional conditions and inequalities and their duals. These are the mixed lattice versions of the lattice properties \eqref{lass} through \eqref{lmod1}, namely the \emph{one-sided associative laws} 

\begin{equation}\label{ass1}
x\lenv (y\lenv z) \leq (x\lenv y)\lenv z, \\[1ex]
\end{equation}
\begin{equation}\label{ass2}
\quad (x\uenv y)\uenv z \leq x\uenv (y\uenv z), \\[2ex]
\end{equation}

the \emph{one-sided distributive laws} 

\begin{equation}\label{dist1}
(z\lenv x)\uenv (z\lenv y) \leq z\lenv(x\uenv y), \\[1ex]
\end{equation}
\begin{equation}\label{dist2}
z\uenv(x\lenv y) \leq (z\uenv x)\lenv (z\uenv y), \\[2ex]
\end{equation}

and the \emph{modular inequalities}

\begin{equation}\label{mod1}
(z\lenv x)\uenv (z\lenv y) \leq z\lenv((z\lenv x)\uenv y), \\[1ex]
\end{equation}
\begin{equation}\label{mod2}
(z\uenv x)\lenv (z\uenv y) \geq z\uenv((z\uenv x)\lenv y). \\[2ex]
\end{equation}

The following dual statements are called the \emph{quasi-regularity conditions}.

\begin{equation}\label{quasireg1}
x\sleq z \; \textrm{ and } \; y\sleq z  \quad \implies \quad x\uenv y \sleq z \\[1ex]
\end{equation}
\begin{equation}\label{quasireg2}
z\sleq x \; \textrm{ and } \; z\sleq y \quad \implies \quad z\sleq x\lenv y. \\[2ex]
\end{equation}

The next two theorems are the main result of this section, and they establish the relationships between 
the preceding conditions. 

\begin{thm}\label{main}
Let $\Mm=(M,\leq,\sleq)$ be a mixed lattice. The following statements are equivalent.
\begin{enumerate}[\normalfont(a)]
\item 
The quasi-regularity condition \eqref{quasireg1} holds in $\Mm$.
\item 
The modular inequality \eqref{mod1} holds in $\Mm$.
\item 
If $x\sleq z$ then $x\uenv (z\lenv y)\leq z\lenv (x\uenv y)$ holds for all $y\in\Mm$. 
%
\item 
The quasi-regularity condition \eqref{quasireg2} holds in $\Mm$.
%
%
\item 
The modular inequality \eqref{mod2} holds in $\Mm$.
\end{enumerate}
Moreover, each of the above conditions implies the pre-regularity condition \eqref{prereg}, and hence the equivalences given in \eqref{r0}.
\end{thm}

\begin{bproof}
(a)$\implies$(b) \, 
First we note that $z\lenv x \sleq z$ and $z\lenv y \sleq z$, so by assumption we have $(z\lenv x)\uenv (z\lenv y)\sleq z$. 
Moreover, $z\lenv y \leq y$ and $z\lenv x \sleq z\lenv x$ imply that $(z\lenv x)\uenv (z\lenv y)\leq (z\lenv x)\uenv y$, by Theorem \ref{basic}(c). Combining these, and using Theorem \ref{basic}(c) again, we obtain the desired inequality  
$(z\lenv x)\uenv (z\lenv y)\leq z\lenv ((z\lenv x)\uenv y).$

(b)$\implies$(c) \, 
If $y\in \Mm$ and $x\sleq z$ then $z\lenv x=x$, and substituting this into \eqref{mod1} gives 
$x\uenv (z\lenv y)\leq z\lenv (x\uenv y)$.

(c)$\implies$(d) \, 
Let $x\sleq z$. If $x\sleq y$ then $x\uenv y=y$, so by hypothesis and Theorem \ref{basic}(b) we have
\[
z\lenv y \leq x\uenv (z\lenv y) \leq z\lenv (x\uenv y) =z\lenv y. 
\]
Hence, $z\lenv y = x \uenv (z\lenv y)$, and by Theorem \ref{basic}(d) this is equivalent to $x\sleq z\lenv y$. This shows that (d) holds. 
%

The implications (d)$\implies$ (e) $\implies$(c) $\implies$(a) are proved similarly, as they are duals of the implications already proved.

To complete the proof, we will show that the quasi-regularity condition (a) implies \eqref{prereg}. For this, assume that (a) holds and let $x\sleq y$. Then, since $x\sleq x$, we have $x \sleq x\lenv y \sleq x$. Hence, by antisymmetry $x=x\lenv y\leq y$, and so \eqref{prereg} holds.  
It was already shown in Theorem \ref{prereg_r0} that the pre-regularity condition implies the equivalences in \eqref{r0}.
\end{bproof}


Next we consider the one-sided distributive and associative laws.

\begin{thm}\label{main1}
Let $\Mm=(M,\leq,\sleq)$ be a mixed lattice. The following statements are equivalent.
\begin{enumerate}[\normalfont(a)]
\item 
If $x\sleq y$ then $z\lenv x \sleq z\lenv y$ holds for all $z\in\Mm$.
\item 
The one-sided distributive law \eqref{dist1} holds in $\Mm$.
\item 
The one-sided associative law \eqref{ass1} holds in $\Mm$.
\end{enumerate}
%
Similarly, the following dual statements are also mutually equivalent.
\begin{enumerate}[\normalfont(d)]
\item 
If $x\sleq y$ then $z\uenv x \sleq z\uenv y$ holds for all $z\in\Mm$.
\item[\normalfont(e)] 
The one-sided distributive law \eqref{dist2} holds in $\Mm$.
\item[\normalfont(f)] 
The one-sided associative law \eqref{ass2} holds in $\Mm$.
\end{enumerate}
Moreover, each of the above conditions implies the equivalent conditions of Theorem \ref{main}.
\end{thm}

\begin{bproof}
%
(a)$\implies$(b) \, 
Since $x\sleq x\uenv y$, we have $z\lenv x \sleq z\lenv (x\uenv y)$ by (a). On the other hand, $z\lenv y\sleq z$ and $z\lenv y \leq y\leq x\uenv y$ imply by Theorem \ref{basic}(c) that $z\lenv y \leq z\lenv (x\uenv y)$. Combining these using Theorem \ref{basic}(c) we obtain %
$(z\lenv x)\uenv (z\lenv y) \leq [z\lenv(x\uenv y)]\uenv [z\lenv(x\uenv y)]=z\lenv(x\uenv y)$.

(b)$\implies$(a) \, 
Let $z\in \Mm$ and $x\sleq y$, or equivalently, $x\uenv y=y$. Then by the distributive inequality we have 
\[
z\lenv y \leq (z\lenv x) \uenv (z\lenv y) \leq z\lenv (x\uenv y) = z\lenv y,
\]
hence $(z\lenv x) \uenv (z\lenv y)=z\lenv y$, by antisymmetry. By Theorem \ref{basic}(d) this is equivalent to $z\lenv x \sleq z\lenv y$.

(a)$\implies$(c) \, 
Since $y\lenv z\sleq y$, we have $x\lenv(y\lenv z)\sleq x\lenv y$ by (a). We also have $x\lenv(y\lenv z)\leq y\lenv z\leq z$, so by Theorem \ref{basic}(c) we obtain $x\lenv (y\lenv z) \leq (x\lenv y)\lenv z$.

(c)$\implies$(a) \, 
Let $z\sleq y$, or equivalently, $y\lenv z=z$. Then, since $x\lenv y\sleq x$, we have by Theorem \ref{basic}(c) that $(x\lenv y)\lenv z \leq x\lenv z$. Using the one-sided associativity law then gives
\[
x\lenv z = x\lenv (y\lenv z) \leq (x\lenv y)\lenv z \leq x\lenv z.
\]
Thus we have $x\lenv z = (x\lenv y)\lenv z \sleq x\lenv y$, and so (a) holds.

The equivalence of (d), (e) and (f) is proved similarly, as they are duals of the statements already considered.

To finish the proof, we show that (c) implies the condition (a) of Theorem \ref{main}. 
Let $z\sleq x$ and $z\sleq y$, or equivalently, $x\lenv z=z$ and $y\lenv z=z$. Then, using the one-sided associativity law we obtain 
\[
z\geq (x\lenv y)\lenv z \geq x\lenv (y\lenv z) = x\lenv z=z.
\]
Thus, $(x\lenv y)\lenv z=z$, or equivalently, $z\sleq x\lenv y$. Hence the quasi-regularity condition \eqref{quasireg1} holds, and the proof is complete.
\end{bproof}


\begin{remark2}\label{huom111}
In algebraic mixed lattice structures, the implications (a)$\implies$(b), (b)$\implies$(c) and (b)$\implies$(d) with their dual versions were proved in \cite[Theorems 3.3 and 3.5]{ars} for mixed lattice semigroups, and in \cite{jj2} for mixed lattice groups. It was also proved in \cite[Theorem 3.2]{jj2} that in mixed lattice groups the quasi-regularity condition \eqref{quasireg1} implies the condition (a) of Theorem \ref{main1}. Therefore, all the statements of Theorems \ref{main} and \ref{main1} are equivalent in mixed lattice groups. 
However, it is not known whether these statements are equivalent in a more general case, and this remains an open question at this point. 
\end{remark2}




Many examples of mixed lattices in which the conditions of Theorems \ref{main} and \ref{main1} hold can be found in \cite{jj1}. 
We now consider a couple of additional examples.

\begin{example2}\label{esim1}
Let $M=\N$ and let $\leq$ be the usual ordering of natural numbers. Define another partial ordering $\sleq$ by $x\sleq y$ if $x$ divides $y$. Then 
$\Mm=(M,\leq,\sleq)$ is a mixed lattice in which the pre-regularity condition \eqref{prereg} clearly holds, but none of the conditions of Theorems \ref{main} and \ref{main1} hold. 
Indeed, if $x=2$, $y=4$ and $z=6$ then $x\sleq y$, $x\sleq z$ 
and $z\lenv y=3$, but $x\sleq z\lenv y$ does not hold. Hence, the condition of Theorem \ref{main}(a) does not hold in $\Mm$. 
In this example, the partially ordered sets $(M, \leq)$ and $(M, \sleq)$ are both distributive lattices, and this shows that the associative, modular and distributive inequalities do not necessarily hold even if $(M,\leq,\sleq)$ is a distributive lattice with respect to both partial orderings (although in this case, the corresponding properties \eqref{lass}--\eqref{lmod1} hold for both partial orderings separately). It should be noted, however, that a mixed lattice is not necessarily a lattice with respect to either partial ordering (see \cite[Example 2.17]{jj1}). 

As with lattices, we can draw diagrams to describe mixed lattices. An example is given in Figure \ref{fig:kuva}, which represents a mixed sublattice $\Ss=\{1,2,3,4,6,12\}$ of $\Mm$, where $1$ is the smallest element and $12$ is the largest element.

\begin{figure}[!h]
\centering
\begin{tikzpicture}[scale=.7]
  \node (one) at (90:3.7cm) {$12$};
  \node (b) at (160:2.2cm) {$4$};
  \node (a) at (230:3.2cm) {$2$};
  \node (zero) at (270:3.6cm) {$1$};
  \node (c) at (345:1.8cm) {$3$};
  \node (d) at (55:3.05cm) {$6$};
  \draw [semithick] (zero) -- (a) -- (b) -- (one) -- (d) -- (c) -- (zero);
	\draw [semithick] (d) -- (a);
	\draw [semithick, dashed] (b) -- (c);
	\draw [semithick, dashed] (b) -- (d);
	\draw [semithick, dashed] (a) -- (c);
\end{tikzpicture}
\caption{\quad }\label{fig:kuva}
\end{figure}

In Figure \ref{fig:kuva}, an ascending solid line from element $x$ to $y$ means $x\leq y$ and $x\sleq y$, while an ascending dashed line means $x\leq y$ (but not $x\sleq y$).
\end{example2}

\blue{
Our next example shows that the mixed lattice of the preceding example can be extended to a mixed lattice group that is closely related to the theory of divisibility. This also provides an example of a mixed lattice group in which the group operation is a multiplication.} 

\begin{example2}\label{multiplicative}
Let $\Mm=(\N,\leq,\sleq)$ be the mixed lattice of the preceding example. If we consider $\Mm$ together with the ordinary multiplication of integers, then $(\Mm, \cdot)$ is a multiplicative semigroup but 
not a mixed lattice semigroup, because in a mixed lattice semigroup the semigroup operation is distributive over the mixed envelopes (i.e. $(x\lenv y)\cdot z=(x\cdot z)\lenv (y\cdot z)$ and $(x\uenv y)\cdot z=(x\cdot z)\uenv (y\cdot z)$) \cite[pp. 16]{ars}, but this does not hold in $(\Mm, \cdot)$. For example, $7\lenv 5=1$ but $(2\cdot 7)\lenv (2\cdot 5)=14\lenv 10=7 \neq 2\cdot 1$. Arsove and Leutwiler \cite{ars} have defined a more general structure, called an \emph{$A$-structure}, 
\blue{which is a partially ordered commutative semigroup $S$ 
with neutral element and two partial orders such that all elements are positive 
and the mixed upper envelope $x\uenv y$ exists in $S$ for all $x,y\in S$ (see also \cite[Theorem 2.2]{eri} for a characterization of $A$-structures).} 
These conditions hold in $(\Mm, \cdot)$, and hence $(\Mm, \cdot)$ is an $A$-structure. 

Now we can apply to $(\Mm, \cdot)$ the well-known construction of positive rational numbers $\mathbb{Q}_+$ by extending $\N$ to a semifield of quotients \blue{(see e.g. \cite[Chapter II.2]{hebi})}. Eriksson-Bique has shown in \cite[Theorem 3.4]{eri} that if $\Aa$ is an $A$-structure, then essentially the same procedure extends $\Aa$ to a mixed lattice group  (except that now we are dealing with just one operation, namely multiplication). It follows that $\Gg=(\mathbb{Q}_+ , \cdot, \leq, \sleq)$ is a mixed lattice group. In the present case, the partial orders $\leq$ and $\sleq$ in the group extension $\Gg$ are defined by 
\[
\frac{a}{b}\leq \frac{c}{d} \, \textrm{ in } \, \Gg  \,\iff \, ad\leq bc \, \textrm{ in } \, (\Mm, \cdot)
\]
and
\[
\frac{a}{b}\sleq \frac{c}{d} \, \textrm{ in } \, \Gg  \,\iff \, ad\sleq bc \, \textrm{ in } \, (\Mm, \cdot).
\] 
\vspace{1mm} \\
This example is related to the theory of divisibility, where divisibility of rational numbers can be defined as follows: if $q,r\in \mathbb{Q}$ then $q$ divides $r$ if $r/q$ is an integer. By the preceding discussion, in $\Gg$, the element $s=r\lenv q$ is the largest rational number such that $s\leq q$ and $s$ divides $r$. Similarly, $t=r\uenv q$ is the smallest rational number such that $t\geq q$ and $r$ divides $t$. As noted above, the identities such as $z\cdot (x\lenv y)=(z\cdot x)\lenv (z\cdot y)$ do not hold in $(\Mm, \cdot)$ but this problem disappears in the group extension. For instance, if we return to the above example and let $r=7\lenv 5$, then we are looking for the largest rational number $r$ such that $r\leq 5$ and $7 = rk$ with $k\in \N$. The largest such $r$ is associated with the smallest possible $k$, and for $k=1$ we have $r=7> 5$, and for $k=2$ we have $r=7/2\leq 5$. Hence, $7\lenv 5=7/2$, and similarly, $14\lenv 10=7$, so  $(2\cdot 7)\lenv (2\cdot 5)=2\cdot (7\lenv 5)$. In the next section, we will show that this distributivity property of the group operation is essential in the axiomatic definition of a mixed lattice group. 
%
\end{example2}

The next theorem gives conditions under which equality holds in the associative and distributive  
inequalities \eqref{ass1} and \eqref{dist1}. Dual statements hold for \eqref{ass2} and \eqref{dist2}.

\begin{thm}\label{equalities}
Let $\Mm=(M,\leq,\sleq)$ be a mixed lattice such that the equivalent conditions of Theorem \ref{main1} hold. Then
\begin{enumerate}[\normalfont(a)]
\item
If $y\sleq x$ or $z\sleq y$ or $y\leq z$ then 
$x\lenv (y\lenv z)=(x\lenv y)\lenv z$.
\item
If $x\sleq y$ or $y\leq x$ then  
$(z\lenv x)\uenv (z\lenv y) = z\lenv(x\uenv y)$. 
\end{enumerate}
\end{thm}

\begin{bproof}
(a) \,
If $y\sleq x$ then $x\lenv y=y$ by Theorem \ref{basic}(d). Similarly, $y\lenv z\sleq y\sleq x$, so $x\lenv (y\lenv z)=y\lenv z$. Hence 
we have $(x\lenv y)\lenv z =y\lenv z = x\lenv (y\lenv z)$.

If $z\sleq y$ then $y\lenv z=z$, and since $x\lenv y\sleq x$, using the one-sided associative law \eqref{ass1} and Theorem \ref{basic}(c) we get
\[
x\lenv z=x\lenv(y\lenv z)\leq (x\lenv y)\lenv z \leq x\lenv z.
\]
This implies that $x\lenv (y\lenv z)=(x\lenv y)\lenv z$.

If $y\leq z$ then $x\lenv y\leq y\leq z$, so by \eqref{r0} we have $y\lenv z=y$ and  
\[
x\lenv(y\lenv z)=x\lenv y =  (x\lenv y)\lenv z. 
\]

(b) \,
If $x\sleq y$ then $x\uenv y=y$ by Theorem \ref{basic}(d). By Theorem \ref{main1}(a) we also have $z\lenv x\sleq z\lenv y$ for all $z\in\Mm$. By Theorem \ref{basic}(d) this is equivalent to $(z\lenv x)\uenv (z\lenv y) = z\lenv y =z\lenv(x\uenv y)$.

If $y\leq x$ then $x\uenv y=x$ by \eqref{r0}, and by Theorem \ref{basic}(c) $z\lenv y\leq z\lenv x$. So by Theorem \ref{basic}(d) we have 
$
z\lenv(x\uenv y) = z\lenv x = (z\lenv x)\uenv (z\lenv y), %
$  
which is the desired result. 
\end{bproof}

\section{Mixed lattices as algebras}
\label{sec:sec3}

In the preceding section we defined a mixed lattice as a partially ordered set. In this section, we give an alternative definition by considering a mixed lattice as an algebra with two binary operations. A similar alternative approach is used in lattice theory, where lattices can be defined as algebras with two binary  operations $\vee$ and $\wedge$ that are commutative, idempotent and associative. In mixed lattices we have two partial orderings, and the mixed lattice operations $\lenv$ and $\uenv$ are neither commutative nor associative, and consequently, our postulates for mixed lattices are a bit more complicated. The task of determining the appropriate postulates is guided by the results of Section \ref{sec:sec2}. 

Let $M$ be a set with two binary operations $\lenv$ and $\uenv$ (not necessarily commutative or associative), and consider the algebra $(M,\uenv,\lenv)$.

\begin{defn}\label{mlalgebra}
\blue{An algebra $\Mm=(M,\uenv,\lenv)$ is called an \emph{algebraic mixed lattice}} if the following conditions hold for all $x,y,z\in M$. \\
\[ 
\begin{array}{lcr}
x\uenv y  =  x\lenv y \quad \implies \quad x=y & & \mathrm{(M1)}  \\ 
(x\lenv y)\uenv x = x \quad \textrm{and} \quad (x\uenv y)\lenv x = x  & & \mathrm{(M2a)} \\
x \uenv (y\lenv x) = x \quad \textrm{and} \quad x \lenv (y\uenv x) = x & & \mathrm{(M2b)} \\
z\lenv (x\uenv y) \; = \; [\,z\lenv (x\uenv y)\,]\uenv (z\lenv y)  & & \mathrm{(M3a)}  \\
z\uenv (x\lenv y) \; = \; [\,z\uenv (x\lenv y)\,]\lenv (z\uenv y)  & & \mathrm{(M3b)}  \\
(x\uenv y)\lenv z \; = \; [\,(x\uenv y)\lenv z\,] \uenv (x\lenv z) & & \mathrm{(M4a)}  \\
(x\lenv y)\uenv z \; = \; [\,(x\lenv y)\uenv z \,] \lenv (x\uenv z) & & \mathrm{(M4b)}  \\
\end{array} 
\]
\end{defn}
\vspace{1mm}

\begin{remark2}
\blue{Regarding Definition \ref{mlalgebra}, it should be noted that from the universal algebra point of view, the conditions (M2a)--(M4b) are identities, whereas (M1) is a quasi-identity. Since identities are also quasi-identities, it follows that the class of mixed lattices as defined here forms a quasi-variety. This immediately raises the following question: Is it possible to replace the condition (M1) by identities? In other words, is the class of mixed lattices a variety?} 

As another remark we note that the postulates (M3a) and (M4a) can be combined to get 
$
(x\lenv v) \lenv [\, (x\uenv y)\lenv (u\uenv v) \,] \, = \, x\lenv v.
$
Similarly, (M3b) and (M4b) can be combined into  
$
(x\uenv v) \uenv [\, (x\lenv y)\uenv (u\lenv v) \,] \, = \, x\uenv v.
$
\end{remark2}

Of course, we would like to show that the above definition is equivalent to Definition \ref{ml} given in Section \ref{sec:sec2}. Let us first derive some easy consequences of the postulates.

\begin{thm}\label{oper_prop}
Let $\Mm=(M,\uenv,\lenv)$ be an algebraic mixed lattice. Then the following hold for all $x,y,z\in M$.
\begin{enumerate}[\normalfont(a)]
\item
$x\lenv x=x$ \, and \, $x\uenv x=x$.
\item
$x=y\lenv x$  $\iff$  $y=x\uenv y$ \quad and \quad $y=y\lenv x$  $\iff$  $x=x\uenv y$.
\item
$x\lenv (x\lenv y) = x\lenv y$ \, and \, $(x\lenv y)\lenv y = x\lenv y$.
\item
$x\uenv (x\uenv y) = x\uenv y$ \, and \, $(x\uenv y)\uenv y = x\uenv y$.
\end{enumerate}
\end{thm}

\begin{bproof}
(a)\; Applying the first identity of (M2b) and the second identity of (M2a) gives 
$x\lenv x = [ x\uenv (y\lenv x)] \lenv x = x$. 
The identity $x\uenv x=x$ is proved similarly.

(b)\; 
If $y\lenv x=x$ then by (M2a) we have $y=(y\lenv x)\uenv y = x\uenv y$, and the reverse implication is proved similarly. The other equivalence is similar too. 

(c)\; 
By the first identity of (M2a) we have $x=(x\lenv y)\uenv x$, and this implies 
\[
x\lenv (x\lenv y) = [(x\lenv y)\uenv x]\lenv (x\lenv y) = x\lenv y,
\]
where we used the second indentity of (M2a). The other identity and the identities in (d) can be proved similarly.  
\end{bproof}

We now show that Definition \ref{mlalgebra} and the definition of a mixed lattice as a partially ordered set together with the property \eqref{r0} are equivalent in the sense of the following theorem.

\begin{thm}\label{main2}
Let $\Mm=(M, \uenv,\lenv)$ be an algebraic mixed lattice and define relations $\sleq$ and $\leq$ by 
\begin{equation}\label{orders2}
x\sleq y \,\iff \, y\lenv x=x \quad \textrm{ and } \quad x\leq y \, \iff \, x\lenv y=x. \\[1ex]
\end{equation}
Then $\sleq$ and $\leq$ are partial orders, and the associated partially ordered set $\Mm_p=(M,\leq,\sleq)$ is a mixed lattice such that \eqref{r0} holds. 

Conversely, if a partially ordered set $\Mm=(M,\leq,\sleq)$ is a mixed lattice such that \eqref{r0} holds then the associated algebra $\Mm_a=(M, \uenv,\lenv)$ is an algebraic mixed lattice, where $\uenv$ and $\lenv$ are the mixed envelopes of $\Mm$.
\end{thm}

\begin{bproof}
By Theorem \ref{oper_prop}(a) we have $x\lenv x=x$, and so $x\leq x$ and $x\sleq x$ holds for all $x$. Moreover, if $x\sleq y$ and $y\sleq x$ then $y=x\lenv y$ and $x=y\lenv x$. But by Theorem \ref{oper_prop}(b) the latter is equivalent to $y=x\uenv y$, and so we have $y=x\lenv y=x\uenv y$. By (M1) this implies that $x=y$. A similar argument shows that if $x\leq y$ and $y\leq x$ then $x=y$. Hence $\sleq$ and $\leq$ are both reflexive and antisymmetric. It remains to show that they are transitive. 
 
If $z\leq y$ and $y\leq x$ then $z\lenv y =z$ and $x\uenv y=x$, by Theorem \ref{oper_prop}(b). Substituting these into (M3a) gives $z\lenv x=(z\lenv x)\uenv z$. By (M2a) the right-hand side equals $z$, and so $z=z\lenv x$, or $z\leq x$. Therefore $\leq$ is transitive, thus proving that $\leq$ is a partial order. 

Assume next that $z\sleq x$ and $x\sleq y$. Then $x\lenv z =z$ and $y\lenv x=x$, or $x\uenv y=y$, by Theorem \ref{oper_prop}(b).  Substituting these into (M4a) yields $y\lenv z = (y\lenv z) \uenv z$. On the other hand, by Theorem \ref{oper_prop}(c) we have $y\lenv z = (y\lenv z) \lenv z$. Hence, $(y\lenv z) \lenv z = (y\lenv z) \uenv z$, and it follows by (M1) that $(y\lenv z) = z$. Thus, $z\sleq y$ and so $\sleq$ is transitive and therefore a partial order. 

We still need to show that the operations $\uenv$ and $\lenv$ are indeed the mixed upper and lower envelopes, as defined by \eqref{upperenv} and \eqref{lowerenv}. Let $A=\{\,w\in M:  w\sleq x \, \textrm{ and } \, w\leq y \,\}$. By (M2a) and (M2b) we have $x=(x\lenv y)\uenv x$ and $y=y\uenv(x\lenv y)$, and so $x\lenv y\sleq x$ and $x\lenv y\leq y$, by the definition of $\leq$ and $\sleq$. This shows that $x\lenv y\in A$. Suppose that $z\sleq x$ and $z\leq y$, that is, $z=x\lenv z$ and $y=y\uenv z$. Then by (M3a) we have 
\[
x\lenv y = x\lenv (y\uenv z) = [\, x\lenv (y\uenv z) \,] \uenv (x\lenv z) = (x\lenv y)\uenv z,
\]
or $z\leq x\lenv y$. Hence $x\lenv y =\max A$, as required. 
Similarly, $x\uenv y =\min \{\,w\in M:  w\sgeq x \, \textrm{ and } \, w\geq y \,\}$.

The converse statement was mostly proved in Section \ref{sec:sec2} (Theorems \ref{basic} and \ref{prereg_r0}) where we showed that the mixed envelopes have all the properties defined by the postulates in Definition \ref{mlalgebra}. First we note that by \eqref{r0} and Theorem \ref{basic}(d), the relationship between the partial orders and the operations $\uenv$ and $\lenv$ are those given in \eqref{orders2}. Then the postulates (M2a) and (M2b) are just the absorbtion laws of Theorems \ref{basic}(e) and \ref{prereg_r0}. 
\blue{For (M3a) we note that $z\sleq z$ and $y\leq x\uenv y$ (by Theorem \ref{basic}(b)), and so by Theorem \ref{basic}(c) we have $z\lenv y \leq z\lenv (x\uenv y)$. By \eqref{r0} this inequality is equivalent to (M3a). Similarly, $x\sleq x\uenv y$ and $z\leq z$, so again by Theorem \ref{basic}(c) it follows that $x\lenv z\leq (x\uenv y)\lenv z$, which is equivalent to (M4a). The dual statements (M3b) and (M4b) hold by similar arguments.} 
As for (M1), we note that if $x\lenv y=x\uenv y$ then by Theorem \ref{basic}(b) we have $y\leq x\uenv y=x\lenv y \leq y$ and $x\sleq x\uenv y=x\lenv y \sleq x$, and it follows by antisymmetry that $x=y$. 
\end{bproof}

An algebraic mixed lattice $\Mm=(M,\uenv,\lenv)$ is called \emph{quasi-regular} if the following additional conditions holds in $\Mm$: \\
\[ 
\begin{array}{lcr}
z\lenv (x\uenv y) \, = \, (z\lenv x) \uenv [\,z\lenv (x\uenv y)\,]  &  & \quad \mathrm{(Q1)}  \\
z\uenv (x\lenv y) \, = \, (z\uenv x) \lenv [\,z\uenv (x\lenv y)\,]  &  & \quad \mathrm{(Q2)}  \\
\end{array} 
\]
\vspace{1mm}

If we use \eqref{orders2} to define two partial orders, then the postulates (Q1) and (Q2) are just restatements of the conditions (a) and (d) of Theorem \ref{main1}. 
By Theorem \ref{main1} these conditions could be stated in many equivalent ways. We have chosen the above formulation because it has a similar form as the postulates (M3a)--(M4b).

The weaker pre-regularity condition \eqref{prereg} can be stated as follows:
\[ 
\begin{array}{lcr}
y\lenv x =x \quad \implies \quad x\lenv y= x & & \quad \mathrm{(P)}  \\
\end{array} 
\]
It is easy to see that (Q1) (or (Q2)) implies (P). Indeed, if (Q1) holds and $y\lenv x =x$ then by Theorem \ref{oper_prop}(b) $x\uenv y=y$, and substituting this into (Q1) and putting $z=x$ we obtain $x\lenv y = x\uenv (x\lenv y)$. By Theorem \ref{oper_prop}(c) we have $x\lenv y=x\lenv (x\lenv y)$, and so $x\uenv (x\lenv y)=x\lenv (x\lenv y)$. By (M1) this implies that $x\lenv y=x$, and so (P) holds. 
If the partial orders $\sleq$ and $\leq$ are defined as in \eqref{orders2} then this conclusion is precisely the same as in Theorem \ref{main}. 

Using the axiomatic definition of a mixed lattice, we can give an alternative description of a mixed lattice group as well. Definition \ref{mlg} defines a mixed lattice group as a partially ordered group that is also a mixed lattice. The following is the algebraic definition.

\begin{defn}\label{mlgroup}
Let $(G,*)$ be a commutative group with two additional binary operations $\uenv$ and $\lenv$ such that the algebra $(G,\uenv,\lenv)$ is an algebraic mixed lattice. Then $(G,*,\uenv,\lenv)$ is called a \emph{mixed lattice group} if the group operation $*$ is distributive over $\uenv$ and $\lenv$ , that is, if %
$$
(x\uenv y) * z = (x*z)\uenv (y*z) \quad \textrm{and} \quad (x\lenv y) * z = (x*z)\lenv (y*z)  
$$
hold for all $x,y,z\in G$.
\end{defn}

The next result shows that the two definitions \ref{mlg} and \ref{mlgroup} of a mixed lattice group are equivalent.  

\begin{thm}
If $\Gg=(G,*,\uenv,\lenv)$ is a mixed lattice group and the partial orders $\leq$ and $\sleq$ are defined by \eqref{orders2} then the partially ordered group $\Gg_p=(G,*,\leq,\sleq)$ is a mixed lattice group such that \eqref{r0} holds.
 
Conversely, if the partially ordered group $\Gg=(G,*,\leq,\sleq)$ is a mixed lattice group such that \eqref{r0} holds, and $\uenv$ and $\lenv$ are the mixed envelopes in $\Gg$ then $\Gg_a=(G,*,\uenv,\lenv)$ is a mixed lattice group.
\end{thm}

\begin{bproof} 
By Theorem \ref{main2} we only need to show that if $\Gg=(G,*,\uenv,\lenv)$ is a mixed lattice group then the partial orders defined by \ref{orders2} are group orderings, that is, they satisfy the condition \eqref{gorder}. To see this, let $x\leq y$ and $z\in G$. Then $x=x\lenv y$ and by Definition \ref{mlgroup} we have $x*z=(x\lenv y)*z=(x*z)\lenv (y*z)$. Hence $x*z\leq y*z$. A similar argument shows that $x\sleq y$ implies $x*z\sleq y*z$. Thus, the partially ordered group $\Gg_p=(G,*,\leq,\sleq)$ is a mixed lattice group. 

Conversely, let the partially ordered group $\Gg=(G,*,\leq,\sleq)$ be a mixed lattice group such that \eqref{r0} holds. It has been shown in \cite[Lemma 3.1]{eri} that for every $x,y,z\in G$ the identities $(x\uenv y) * z = (x*z)\uenv (y*z)$ and $(x\lenv y) * z = (x*z)\lenv (y*z)$ hold, and so it follows by Theorem \ref{main2} that $\Gg_a=(G,*,\uenv,\lenv)$ is a mixed lattice group.
\end{bproof}





\bibliographystyle{plain}

\end{document}